\documentclass[12pt]{article}
\usepackage{amssymb,amsbsy,amsmath,amsfonts,amssymb}
\usepackage{latexsym,euscript,exscale}

\usepackage{times}

\newcommand{\pf}{

\smallskip

\noindent {\it Proof : }}

\newcommand {\N}{\mathbb N}

\newcommand {\Q}{\mathbb Q}

\newcommand{\pff}{$\hfill \square$
\smallskip}

\newcommand{\norm}[1]{\ensuremath{\left\|#1\right\|}}

\newtheorem{prop}{Proposition}[section]
\newtheorem{lemm}[prop]{Lemma}
\newtheorem{theo}[prop]{Theorem}
\newtheorem{conj}[prop]{Question}
\newtheorem{coro}[prop]{Corollary}
\newtheorem{rema}[prop]{Remark}
\newtheorem{defi}[prop]{Definition}

\title{Minimality, homogeneity and topological 0-1 laws for
subspaces of a Banach
space}

\author{Valentin Ferenczi \footnote{Part of this paper was written at the University of
S\~ao Paulo under FAPESP Grant 02/09662-1.}}

\date{ }

\begin{document}

\maketitle

\begin{abstract}
If a Banach space is saturated with basic sequences whose linear span
embeds into the linear span
of any subsequence, then it contains a minimal subspace. It follows
that any Banach space is either
ergodic or contains a minimal subspace.

 For a Banach space $X$ with
an (unconditional) basis,
topological 0-1 law type dichotomies are stated for block-subspaces
of $X$ as well as for subspaces of
$X$ with a successive FDD on its basis. A uniformity principle for
properties of block-sequences,
results about block-homogeneity, and a possible method to construct a
Banach space with an
unconditional basis, which has a complemented subspace without an
unconditional basis, are deduced.
\end{abstract}

\

The starting point of this article is the solution to the
Homogeneous Banach Space Problem given by W.T. Gowers \cite{G} and R.
Komorowski -
N. Tomczak- Jaegermann \cite{KT}. A Banach space is said to be
homogeneous if it is
isomorphic to its infinite dimensional closed subspaces; these authors proved
that a homogeneous Banach space must be isomorphic to $l_2$.

\

Gowers proved that any Banach space with a basis must either have a subspace with an
unconditional basis or a hereditarily indecomposable subspace. By properties of
hereditarily indecomposable Banach spaces, it follows that a homogeneous Banach space must
have an unconditional basis (see e.g. \cite{G} for details about this). Komorowski and
Tomczak-Jaegermann proved that a Banach space with an unconditional basis must contain a
copy of $l_2$ or a subspace with a successive finite-dimensional decomposition on the basis
($2$-dimensional if the space
has finite cotype) which does not have an unconditional basis.
It follows that a homogeneous Banach space must be isomorphic to $l_2$.

While Gowers' dichotomy theorem is based on a general Ramsey type theorem for
block-subspaces in a Banach space with a Schauder basis, the subspace with a
finite-dimensional decomposition constructed in Komorowski and Tomczak-Jaeger\-mann's
theorem can never be isomorphic to a block-subspace. 
If one restricts one's attention to block-subspaces, the standard homogeneous
examples become the sequence spaces $c_0$ and $l_p, 1 \leq p<+\infty$, with their canonical
bases; these spaces are well-known to be isomorphic to their block-subspaces.
Furthermore there are classical theorems which characterize $c_0$ and $l_p, 1 \leq
p<+\infty$ only by means of their block-subspaces. An instance of this is Zippin's theorem
(\cite{LT} Theorem 2.a.9): a normalized basic sequence is perfectly homogeneous (i.e.
equivalent to all its normalized block-sequences) if and only if it is equivalent to the
canonical basis of
$c_0$ or some
$l_p$.  See also
\cite{LT} Theorem 2.a.10.

 So it is very natural
  to ask
  what can be said on the subject of (isomorphic)
 homogeneity restricted to block-subspaces of a given Banach space with a Schauder
basis: if a Banach space $X$ with a basis
$(e_n)_{n \in \N}$  is
isomorphic to its block-subspaces, does it follow that $X$ is isomorphic to
$c_0$ or $l_p, 1 \leq p <+\infty$?
Note that such a basis is not necessarily equivalent to the canonical basis of $c_0$ or
some $l_p$, take $l_2$ with a conditional basis.

\

In the other direction, if a Banach space is not homogeneous, then how many
non-isomorphic subspaces must it contain? This question may be asked in the
setting of the classification of analytic equivalence relations on Polish
spaces by Borel reducibility. This area of research originated from the works of H.
Friedman and L. Stanley \cite{FS} and independently from the works of L. A. Harrington, A.
S. Kechris and A. Louveau \cite{HKL}, and may be thought of as an extension of the notion
of cardinality in terms of complexity, when one compares equivalence relations.

 If
$R$ (resp.
$S$) is an
equivalence relation on a Polish space $E$ (resp. $F$), then it is said
that $(E,R)$ is Borel
reducible to $(F,S)$ if there exists a Borel map $f:E \rightarrow F$ such that
$\forall x,y \in E, xRy \Leftrightarrow f(x)Sf(y)$.
An important equivalence relation is the relation $E_0$: it is
defined on  $2^{\omega}$ by
$$\alpha E_0 \beta \Leftrightarrow \exists m \in \N \forall n \geq m,
\alpha(n)=\beta(n).$$

The relation $E_0$ is a Borel equivalence relation with $2^{\omega}$ classes and which,
furthermore, is
not Borel reducible to equality on $2^{\omega}$, that is, there is no Borel map $f$
from
$2^\omega$ into $2^\omega$ (equivalently, into a Polish space), such that $\alpha E_0 \beta
\Leftrightarrow f(\alpha)=f(\beta)$; such a relation is said to be {\em non-smooth}.
In fact $E_0$ is the $\leq_B$ minimal  non-smooth Borel equivalence relation \cite{HKL}.

 There is a natural way to equip the set of subspaces
of a Banach space
$X$ with a Borel structure \cite{B}, and the relation of isomorphism is analytic in this
setting. The relation
$E_0$ appears to be a natural threshold for results about the relation of isomorphism
between separable Banach spaces
\cite{FG},\cite{FR1},\cite{FR2},\cite{R}.
A Banach space $X$ is said to be {\em ergodic} if $E_0$ is Borel reducible to
isomorphism between subspaces of $X$; in particular, an ergodic Banach space
has continuum many non-isomorphic subspaces, and isomorphism between its
subspaces is non-smooth. The
results in
\cite{B},\cite{FG},\cite{FR1},\cite{FR2},\cite{R} suggest that every Banach
space non-isomorphic to $l_2$ should be ergodic, and we also refer to these
articles for an introduction to the classification of analytic
equivalence relations on Polish spaces by Borel reducibility, and more specifically
to the complexity of isomorphism between Banach spaces.

Restricting our attention to block-subspaces, the
natural question  becomes the following:
if $X$ is a Banach space with a Schauder basis, is it true that either
  $X$ is isomorphic to its
block-subspaces or  $E_0$ is Borel reducible to isomorphism between the
block-subspaces of $X$?

\

Let us provide some ground for this conjecture by noting that, if
we replace isomorphism by equivalence of the
corresponding basic sequences, it is completely solved
by the positive by a result of \cite{FR1} using the theorem of
Zippin: if
$X$ is a Banach space with a normalized basis $(e_n)_{n \in \N}$, then either
$(e_n)_{n \in \N}$ is equivalent to the canonical basis of $c_0$ or
$l_p, 1 \leq p
<+\infty$, or $E_0$ is Borel reducible to equivalence between normalized
block-sequences of $X$.

\

This article is divided in three sections. The results and methods in the first two
sections are mainly independant, although some notation defined in the first section might
be used in the second section. We obtain partial answers to the above conjectures in
various directions; our methods also provide results of combinatorial nature
which are of independant interest in Banach space theory. The third section contains a
refined version of a principle proved in the third section, with an application to the
study of complemented subspaces of a Banach space with an unconditional basis.

\

\

A.M. Pelczar has proved that a Banach space which is
saturated with subsymmetric sequences
contains a minimal subspace \cite{P}.
Our main theorem in this article (Theorem \ref{anyaisomorphic}), proved in
the first section, is the following. If a Banach space $X$ is saturated
with basic sequences whose linear span embed in the
linear span of any subsequence, then $X$ contains a minimal subspace.
In particular, define a basic sequence to be
{\em isomorphically homogeneous} if all subspaces spanned by subsequences
are isomorphic; our result implies that a
Banach space saturated with isomorphically homogeneous basic
sequences contains a minimal subspace.
This is the isomorphic counterpart of Pelczar's result.

In combination with a result of C. Rosendal \cite{R}, it follows
that if $X$ is a Banach space with a Schauder basis, then either $E_0$ is Borel
reducible to isomorphism between block-subspaces of $X$, or $X$ contains
a block-subspace which is block-minimal (i.e. embeds as a block-subspace of
any of its block-subspaces), Corollary
\ref{ergodicorminimal}. This improves a
result of
\cite{FR2} which states that a Banach space contains continuum many
non-isomorphic subspaces or a
minimal subspace.

\

The {\em second topological 0-1 law} (Theorem 8.47 in \cite{ke})
states that in a infinite
product space of Polish spaces, a set with Baire Property
  which is a tail set (i.e. invariant by change of a finite
number of coordinates), is either meager or comeager.
In the second section,  we study the
set $bb_d(X)$ of
  ``rational normalized block-sequences'' of a Banach space $X$ with a Schauder
  basis,
  and a characterization of comeager sets in the natural topology on
$bb_d(X)$ that was
obtained in \cite{FR2}, to deduce a
principle of topological 0-1 law for block-subspaces in $bb_d(X)$
(Theorem \ref{topologicallaw},
Theorem \ref{dichotomytheorem}).

We deduce from this principle a uniformity theorem (remark after
Proposition \ref{theo3}),
and an application to some problems related to the block-homogeneity
property (Proposition \ref{isomorphictolp}).

\

In the third section, we prove a principle of 0-1 topological law in
a Banach space
$X$ with a Schauder basis, for subspaces with a successive finite dimensional
decomposition on the basis
(Proposition \ref{prop20}), again continuing on some work from \cite{FR2}.
  We derive a possible application to a long-standing open question in
Banach space theory: does a
complemented subspace of a Banach space with an unconditional basis
necessarily have an
unconditional basis (Corollary \ref{contrexemple})?

\

Let us fix some notation.
Let $X$ be a Banach space with a Schauder basis $(e_n)_{n \in \N}$.
If $(x_n)_{n \in J}$ is a finite or infinite block-sequence of $X$
then $[x_n]_{n \in \N}$
will stand for its closed linear span.
We shall also use some standard notation about finitely supported
vectors on $(e_n)_{n \in \N}$, for
example, we shall write $x<y$ and say that $x$ and $y$ are successive
when $\max(supp(x)) <
\min(supp(y))$. The set of normalized block-sequences in $X$ is denoted
$bb(X)$. It is a Polish space when equipped with the product of the
norm topology on $X$.

  Let  ${\Q}(X)$ be the set of non-zero
blocks of the basis (i.e. finitely supported vectors) which have
rational coordinates on $(e_n)_{n \in \N}$ (or coordinates in
$\Q+i\Q$ if we deal with a complex
Banach space).
   We
denote by $bb_{\Q}(X)$ the set of block-bases of vectors in ${\Q}(X)$, and by
${\cal G}_{\Q}(X)$
the corresponding set of block-subspaces of $X$.

The notation $bb_{\Q}^{<\omega}(X)$ (resp. $bb_{\Q}^n(X)$) will be used for
the set of finite (resp. length $n$) block-sequences with vectors in
${\Q}(X)$; the set of
finite
block-subspaces generated by a block-sequence in
$bb_{\Q}^{<\omega}(X)$ will be denoted
by $Fin_{\Q}(X)$.

We shall consider $bb_{\Q}(X)$ as a topological space,
when equipped with the product 
of the discrete topology on ${\Q}(X)$.  As ${\Q}(X)$ is countable, this
turns $bb_{\Q}(X)$ into a Polish space. Likewise, ${\Q}(X)^{\omega}$
is a Polish space.

For a finite block
sequence
$\tilde{x}=(x_1,\ldots,x_n) \in bb_{\Q}^{<\omega}(X)$, we denote by
$N_{\Q}(\tilde{x})$
the set of elements of $bb_{\Q}(X)$ whose first $n$ vectors are
$(x_1,\ldots,x_n)$; this is the basic open set associated to $\tilde{x}$.

The set $[\omega]^{\omega}$ is the set of increasing sequences of integers,
which we sometimes identify with infinite subsets
of
$\omega$. It is equipped with
the product of the discrete
topology on $\omega$. The set $[\omega]^{<\omega}$ is the set of
finite increasing sequences of integers. If $a
=(a_1,\ldots,a_k) \in [\omega]^{<\omega}$, then   $[a]$ stands for
the basic open set associated to
$a$, that is the set of increasing sequences of integers of the form
$\{a_1,\ldots,a_k,a_{k+1},\ldots\}$.
If $A \in [\omega]^{\omega}$, then $[A]^{\omega}$ is the set of
increasing sequences of integers in $A$ (where $A$ is seen as a subset of $\omega$).

We recall that two basic sequences $(x_n)_{n \in \N}$ and
$(y_n)_{n \in \N}$ are said to be
{\em equivalent} if the map $T:[x_n]_{n \in \N} \rightarrow [y_n]_{n
\in \N}$ defined by
$T(x_n)=y_n$ for all $n \in \N$ is an isomorphism.
For $C \geq 1$, they are {\em $C$-equivalent} if
$\norm{T}\norm{T^{-1}} \leq C$. A basic sequence is
said to be {\em ($C$-)subsymmetric} if it is ($C$-)equivalent to all its
subsequences.

We shall sometimes use "standard perturbation arguments" without
expliciting them. This
expression will refer to one of the following well-known facts about
block-subspaces of a Banach space $X$
with a Schauder basis. Any basic sequence (resp. block-basic sequence) in $X$
is an arbitrarily small
perturbation of a basic sequence in ${\Q}(X)^{\omega}$ (resp.
block-basic sequence in $bb_{\Q}(X)$),
and in particular is
$1+\epsilon$-equivalent to it, for arbitrarily small $\epsilon>0$.
Any subspace of
$X$ has a subspace which is an arbitrarily small perturbation of a
block-subspace of
$X$ (and in particular, with $1+\epsilon$-equivalence of the
corresponding bases, for arbitrarily
small
$\epsilon>0$). If
$X$ is  reflexive, then any basic sequence in $X$ has
a subsequence which is a perturbation of a block-sequence of $X$ (and
in particular, is
$1+\epsilon$-equivalent to it, for arbitrarily small $\epsilon>0$).

   We shall also use the fact that any Banach space
contains a basic sequence.

\section{Minimal subspaces, isomorphically homogeneous sequences, and
reductions of $E_0$.}

We recall different notions of minimality for Banach spaces.
A Banach space $X$ is said to be {\em minimal} if it embeds into any of
its subspaces. If
$X$ has a basis $(e_n)_{n \in \N}$, then it is said to be {\em
block-minimal} if every block-subspace
of
$X$ has a further block-subspace which is isomorphic to $X$, and
  is said to be
{\em equivalence block-minimal} if every block-sequence of $(x_n)_{n
\in \N}$ has a further
block-sequence which is equivalent to
$(x_n)_{n \in \N}$.

The theorem of Pelczar \cite{P} states that a Banach space
which is saturated with subsymmetric sequences must contain an equivalence
block-minimal subspace with a basis. In this section we prove a version of
her theorem for isomorphism (Theorem \ref{anyaisomorphic}).

We recall that a basic sequence $(x_n)_{n \in \N}$ in a Banach space
is said to be
{\em isomorphically homogeneous} if all subspaces spanned by
subsequences of $(x_n)_{n \in \N}$ are
isomorphic.
  The relevant property for our theorem will be the following
property, which is
obviously weaker than
   the property of being isomorphically homogeneous:   say that a basic
sequence {\em embeds
   (resp. $C$-embeds) into
   its subsequences} if its linear span embeds (resp. $C$-embeds) into the linear
span of any of its subsequences.

\begin{theo}\label{anyaisomorphic}
  A Banach space which is saturated with basic sequences which
   embed into their subsequences contains a minimal subspace.
\end{theo}

For $N \in \N$ let $d_c(N)$ denote an integer such that if $X$ is a
Banach space with a basis
$(e_n)_{n \in \N}$ with basis constant $c$, and $(x_n)_{n \in \N}$
and $(y_n)_{n \in \N}$ are
normalized block-basic sequences of $X$ such that $x_n=y_n$ for all
$n>N$, then $(x_n)_{n \in \N}$
and $(y_n)_{n
\in \N}$ are $d_c(N)$-equivalent.
We leave as an exercise to the reader to check that such an integer exists.

\begin{lemm}\label{firstdiagonalization}  Let $(x_n)_{n \in \N}$ be a
basic sequence in a Banach space
   which embeds into its subsequences. Then there exists $C \geq 1$ and a
   subsequence of $(x_n)_{n \in \N}$ which $C$-embeds into its subsequences.
\end{lemm}

\pf
  Let $(x_n)_{n \in \N}$ be a
basic sequence which embeds into
  its subsequences, and let $c$ be its basis constant.
It is clearly enough to find a subsequence $(y_n)_{n \in \N}$ of
$(x_n)_{n \in \N}$ and $C \geq 1$
such that
$(x_n)_{n \in \N}$ $C$-embeds into any subsequence of $(y_n)_{n \in
\N}$ (with the obvious
definition).

Assuming the conclusion is false, we construct by induction a
sequence of subsequences
$(x_n^k)_{n \in \N}$ of $(x_n)_{n \in \N}$,
  such that
for all $k \in \N$, $(x_n^k)_{n \in \N}$ is a subsequence of
$(x_n^{k-1})_{n \in \N}$ such
that $(x_n)_{n
\in
\N}$ does not
$kd_c(k)$-embeds into
$(x_n^k)_{n
\in
\N}$.

Let $(y_n)_{n \in \N}$ be the diagonal subsequence of $(x_n)_{n \in
\N}$ defined by
$y_n=x_n^n$.
Then $(x_n)_{n \in \N}$ does not $kd_c(k)$-embeds into
$(x^k_1,\ldots,x_{k-1}^k,y_k,y_{k+1},\ldots)$.
  So $(x_n)_{n \in \N}$  does not $k$-embeds in
 $(y_n)_{n \in \N}$. Now
$k$ was arbitrary, so this contradicts our hypothesis. \pff

\begin{lemm}\label{seconddiagonalization}
Let $X$ be a Banach space which is saturated with basic sequences which embed
in their subsequences. Then there exists a subspace $Y$ of $X$ with a
Schauder basis, and a constant
$C \geq 1$ such that every block-sequence of $Y$ (resp. in
$bb_{\Q}(Y)$) has a further
block-sequence (resp. in $bb_{\Q}(Y)$)
which $C$-embeds into its subsequences.
\end{lemm}

\pf  By properties of
hereditarily indecomposable Banach spaces \cite{GM}, a basic sequence
which embeds into its
subsequences cannot span a hereditarily indecomposable space. Thus
$X$ does not contain a hereditarily indecomposable
subspace and by Gowers' dichotomy theorem, we may assume $X$ has an
unconditional basis (let $c$ be
its basis constant).
  If $c_0$ or $l_1$ embeds into $X$ then we are done, so by the
classical theorem
of James, we
  may assume
$X$ is reflexive. Thus by standard perturbation arguments, every
normalized block-sequence in $X$ has
a further normalized block-sequence in
$X$ which embeds into its subsequences (here we also used the obvious
fact that if a basic
sequence $(x_n)_{n \in \N}$ embeds into its subsequences, then so does
any subsequence of $(x_n)_{n \in \N}$).

Assuming the conclusion is false, we construct by induction a
sequence of block-sequences
$(x_n^k)_{n \in \N}$ of $(x_n)_{n \in \N}$, for $k \in \N$,
  such that
for all $k \in \N$, $(x_n^k)_{n \in \N}$
is a block-sequence of $(x_n^{k-1})_{n \in \N}$ such that no block-sequence of
$(x_n^k)_{n \in \N}$  $kd_c(k)^2$-embeds into its subsequences.

Let $(y_n)_{n \in \N}$ be the diagonal block-sequence of $(x_n)_{n
\in \N}$ defined by
$y_n=x_n^n$, and let $(z_n)_{n \in \N}$ be an arbitrary
block-sequence of $(y_n)_{n \in \N}$.

Then
$(x^k_1,\ldots,x_{k-1}^k,z_k,z_{k+1},\ldots)$ is a block-sequence of
$(x_n^k)_{n \in \N}$ and so,
does not $kd_c(k)^2$-embed into its subsequences .
  So $(z_n)_{n \in \N}$  does not $k$-embeds into its subsequences -
this is true as well of its
subsequences. As $k$ was arbitrary, we deduce from Lemma
\ref{firstdiagonalization} that $(z_n)_{n \in
\N}$ does not embed into its subsequences. As
$(z_n)_{n
\in
\N}$ was an arbitrary block-sequence of $(y_n)_{n \in \N}$, this
contradicts our hypothesis.

By standard perturbation arguments, we deduce from this the stated
result with block-sequences in
$bb_{\Q}(Y)$.
\pff

\

Recall that $\Q(X)^{\omega}$ is equipped with the product of the
discrete topology on $\Q(X)$ which
turns it into a Polish space.

\begin{defi} \label{continuouslyembeds} Let $X$ be a Banach space
with a Schauder basis, and let
$(x_n)_{n \in \N} \in {\Q}(X)^{\omega}$. We
   shall say that $(x_n)_{n \in \N}$ {\em continuously embeds
   (resp. $C$- continuously embeds) into its subsequences} if there
exists a continuous map
  $\phi:[\omega]^{\omega} \rightarrow {\Q}(X)^{\omega}$ such for all
$A \in [\omega]^{\omega}$, $\phi(A)$ is a sequence of vectors in
  $[x_n]_{n \in A} \cap {\Q}(X)$ which is equivalent (resp. $C$-equivalent) to
  $(x_n)_{n \in \N}$.
\end{defi}

This definition depends on the Banach space $X$ in
which we pick the basic sequence $(x_n)_{n \in \N}$; this will not cause us
any problem, as it will always be clear which is the underlying space $X$.

\

The
interest of this notion stems from the following lemma, which was
essentially obtained by Rosendal as
part of the proof of \cite{R}, Theorem 11. To prove it, we shall need the following
fact, which is well-known to descriptive set theoricians. The algebra
$\sigma(\Sigma_1^1)$ is the $\sigma$-algebra generated by analytic sets.
 For any $\sigma(\Sigma_1^1)$-measurable function from
$[\omega]^{\omega}$ into a metric space, there exists $B \in [\omega]^{\omega}$ such that
the restriction of $f$ to $[B]^{\omega}$ is continuous.

Indeed, by Silver's Theorem \cite{ke} 21.9,
any analytic set in $[\omega]^{\omega}$ is completely Ramsey, and so any 
$\sigma(\Sigma_1^1)$ set in $[\omega]^{\omega}$ is (completely) Ramsey as well (use for
example
\cite{ke} Theorem 19.14). One concludes using the proof of \cite{MW} Theorem 9.10 which
only uses the Ramsey-measurability of the function.

\begin{lemm}\label{continuoussubsequence} Let $X$ be a Banach space
with a Schauder basis,
let $(x_n)_{n \in \N} \in bb_{\Q}(X)$ be a block-sequence which
$C$-embeds into its
subsequences, and let $\epsilon$ be positive. Then some subsequence
of $(x_n)_{n \in \N}$
$C+\epsilon$-continuously embeds into its subsequences.
\end{lemm}

\pf
By standard perturbation arguments, we may find for each $A \in
[\omega]^{\omega}$ a sequence
$(y_n)_{n \in \N} \in {\Q}(X)^{\omega}$ such that $y_n \in
   [x_k]_{k \in A}$ for all $n \in \N$, and
such that the basic sequences $(x_n)_{n \in \N}$ and $(y_n)_{n \in \N}$ are
$C+\epsilon$-equivalent. The set $P \subset [\omega]^{\omega} \times
{\Q}(X)^{\omega}$ of couples
$(A,(y_n))$ with this property is Borel (even closed), so by the
Jankov - von Neumann
Uniformization Theorem (Theorem 18.1 in \cite{ke}), there exists a
$C$-measurable selector
$f:[\omega]^{\omega} \rightarrow {\Q}(X)^{\omega}$ for $P$. By
the fact before this lemma, there exists
$B \in [\omega]^{\omega}$ such that the restriction of $f$ to
$[B]^{\omega}$ is continuous.
Write $B=(b_k)_{k \in \N}$ where $(b_k)_k$ is increasing. By
composing $f$ with the obviously
continuous maps
$\psi_B:[\omega]^{\omega}
\rightarrow [B]^{\omega}$, defined by $\psi_B((n_k)_{k \in
\N})=(b_{n_k})_{k \in \N}$, and
$\mu_B:{\Q}(X)^{\omega} \rightarrow {\Q}(X)^{\omega}$, defined by
$\mu_B((y_n)_{n \in \N})=(y_{b_n})_{n \in
\N}$, we obtain a continuous map $\phi:[\omega]^{\omega} \rightarrow
{\Q}(X)^{\omega}$ which indicates that $(x_n)_{n \in B}$
$C+\epsilon$-continuously embeds into its
subsequences.
\pff

\

We now start the proof of Theorem \ref{anyaisomorphic}. So we
consider a Banach space $X$ which is
saturated with basic sequences which embed into their subsequences and wish
to find a minimal subspace in $X$.

By Lemma \ref{seconddiagonalization} and Lemma
\ref{continuoussubsequence}, we may assume that $X$ is
a Banach space with a Schauder basis and that there exists $C \geq 1$ such
that every block-sequence in
$bb_{\Q}(X)$ has a further block-sequence in $bb_{\Q}(X)$ which
$C$-continuously embeds into its
subsequences.

  For the rest of the proof $X$ and $C \geq 1$ are fixed with this property.
Recall that the set of block-subspaces of $X$ which are generated by
block-sequences in $bb_{\Q}(X)$
is denoted by ${\cal G}_{\Q}(X)$; the set of finite block-subspaces
which are generated by
block-sequences in $bb_{\Q}^{<\omega}(X)$ is denoted by $Fin_{\Q}(X)$.
If $n \in \N$ and $F \in Fin_{\Q}(X)$ we write $n \leq F$ to mean
that $n \leq \min(supp(x))$ for all $x \in F$.
\

We first express the notion of continuous embedding in terms of a game. 
For $L=[l_n]_{n \in
\N}$ with
$(l_n)_{n \in \N} \in bb_{\Q}(X)$, we define a game
  $H_L$ as follows. A $k$-th move for Player 1 is some 
$n_k \in \N$. A $k$-th
move for Player 2 is some $(F_k,y_k) \in Fin_{\Q}(X) \times {\Q}(X)$, with $n_k \leq F_k
\subset L$ and
  $y_k \in \Sigma_{j=1}^k F_j$.

Player 2 wins the game $H_L$ if  $(y_n)_{n \in \N}$ is $C$-equivalent to $(x_n)_{n \in \N}$.

\

We claim the following:

\begin{lemm}\label{auxiliarygame}
Let $X$ be a Banach space with a Schauder basis, and $(l_n)_{n \in \N} \in
bb_{\Q}(X)$ be a
block-sequence which
$C$-continuously embeds into its subsequences. Let $L=[l_n]_{n \in \N}$.
 Then Player 2 has a
winning strategy in the game $H_L$.
\end{lemm}

\pf Let $\phi$ be the continuous map in Definition \ref{continuouslyembeds}. We describe
a winning strategy for Player 2 by induction.

We assume that Player 1's moves were
$(n_i)_{i \leq k-1}$ and that the $k-1$ first moves prescribed by the winning strategy
for Player
2 were
$(F_i,y_i)_{i
\leq k-1}$, with $F_i$ of the form $[l_{n_i},\ldots,l_{m_i}]$, $n_i \leq m_i$, for all $i
\leq k-1$; letting
$a_{k-1}=[n_1,m_1] \cup \ldots \cup [n_{k-1},m_{k-1}] \in
[\omega]^{<\omega}$, we also assume that
$\phi([a_{k-1}]) \subset N_{\Q}(y_1,\ldots,y_{k-1})$. We now describe the
$k$-th move of the winning strategy for Player 2.

Let $n_k$ be a $k$-th move for Player 1. We may clearly assume that $n_k>m_{k-1}$. Let
$A_k=\cup_{i
\leq k-1}[n_i,m_i])
\cup [n_k,+\infty)
\in [\omega]^{\omega}$. The sequence $\phi(A_k)$ is of the form
$(y_1,\ldots,y_{k-1},y_k,z_{k+1},\ldots)$ for some $y_k, z_{k+1}, \ldots$ in $\Q(X)$.
By continuity of $\phi$ in $A_k$ there exists $m_k>n_k$ such that, if
$a_k=[n_1,m_1] \cup \ldots \cup [n_{k},m_{k}] \in [\omega]^{<\omega}$, then
$\phi([a_{k}]) \subset N_{\Q}(y_1,\ldots,y_{k})$. We may assume that
$\max(supp(x_{m_k})) \geq
\max(supp(y_k))$; so as $y_k \in [x_i]_{i \in A}$, we have that $y_k
\in \oplus_{j=1}^k [x_i]_{i \in [n_j,m_j]}$. So 
$(F_k,y_k)=([l_{n_k},\ldots,l_{m_k}],y_k)$ is an admissible $k$-th move for Player 2 for
which
the induction hypotheses are satisfied.

Repeating this by induction we obtain a sequence $(y_n)_{n \in \N}$
which is equal to $\phi(A)$,
where $A=\cup_{k \in \N}[n_k,m_k]$,
and so which is, in particular, $C$-equivalent to $(x_n)_{n \in \N}$.
\pff

\begin{defi}\label{glm}
Given $L,M$ two block-subspaces in ${\cal G}_{\Q}(X)$,
define the game $G_{L,M}$ as follows. A $k$-th move for Player 1 is some
$(x_k,n_k) \in {\Q}(X) \times \N$,
with $x_k \in L$, and $x_k>x_{k-1}$ if $k \geq 2$. A $k$-th move for
Player 2 is some
$(F_k,y_k) \in Fin_{\Q}(X) \times {\Q}(X)$ with $n_k \leq F_k \subset M$
and $y_k \in F_1 \oplus \ldots \oplus F_k$ for all $k \in \N$.

\[\begin{array}{cccccccc}
                   & I\ :& x_1,n_1 &      &  x_2,n_2 &     & \ldots &    \\
           G_{L,M}\ \ &   &           &      &            &     &        &    \\
                   & II:\ &          & F_1,y_1 &       & F_2,y_2 &   &
\ldots \end{array}\]

Player 2 wins $G_{L,M}$ if
  $(y_n)_{n \in \N}$ is $C$-equivalent to $(x_n)_{n \in \N}$.
\end{defi}

The following easy fact will be needed in the next lemma: if $(x_n)_{n \in \N}$ and
$(y_n)_{n \in \N}$ are $C$-equivalent basic sequences, then for any 
scalar sequence $(\lambda_i)_{i \in \N}$ and sequence $(I_n)_{n \in \N}$ of successive
subsets of
$\N$ such that $\{i \in I_n: \lambda_i \neq 0\} \neq \emptyset, \forall n \in \N$, the basic
sequences
$(\sum_{i
\in I_n}\lambda_i x_i)_{n
\in
\N}$ and
$(\sum_{i \in I_n}\lambda_i y_i)_{n \in \N}$ are $C$-equivalent as well.

\begin{lemm}\label{firststepoftheinduction}
Assume $(l_n)_{n \in \N}$ is a block-sequence  in $bb_{\Q}(X)$ which
$C$-continously embeds into its
subsequences, and let
$L=[l_n, n \in \N]$. Then Player 2 has a winning strategy in the game
$G_{L,L}$.
\end{lemm}

\pf Assume Player 1' first move was
$(x_1,n_1)$;
write $x_1=\sum_{j \leq k_1} \lambda_j l_j$.
Letting in the game $H_L$ Player $1$ play the integer
$n_1$, $k_1$ times,
the winning strategy of Lemma \ref{auxiliarygame} provides moves
$(F_1^1,z_1),\ldots,(F^{k_1}_1,z_{k_1})$
for Player 2 in that game. We let $y_1=\sum_{j \leq k_1} \lambda_j
z_j$, and $F_1=\sum_{j=1}^{k_1} F_1^j$. In particular, $n_1 \leq F_1 \subset L$
and $y_1 \in F_1$.

We describe the choice of $F_p$ and $y_p$ at step $p$.
Assuming Player 1' $p$-th move was
$(x_p,n_p)$;
we write $x_p=\sum_{k_{p-1}< j \leq k_p} \lambda_j l_j$.
Letting in the game $H_L$ Player $1$ play $k_p-k_{p-1}$
times the integer $n_p$,
the winning strategy of Lemma \ref{auxiliarygame} provides moves
$(F_p^{k_{p-1}+1},z_{k_{p-1}+1}),\ldots,(F_p^{k_p},z_{k_p})$ for Player 2
in that game. We let
$y_p=\sum_{k_{p-1}< j \leq k_p} \lambda_j z_j$, and
 $F_p=\sum_{k_{p-1}<j\leq k_p} F_p^j$. In particular, $n_p \leq F_p
\subset L$  and $y_p \in \sum_{j=1}^p F_j$.

Finally by construction,
$(z_n)_{n \in \N}$ is $C$-equivalent to $(l_n)_{n \in \N}$. It follows that
$(y_p)_{p \in \N}$ is $C$-equivalent to $(x_p)_{p \in \N}$. \pff

\

The non-trivial Lemma \ref{firststepoftheinduction} will serve as the
first step of a final
induction which is on the model of the demonstration of Pelczar in
\cite{P} (note
that there, the first step of the induction was straightforward). The
rest of our reasoning in this
section will now be along the lines of her work, with the difference
that we chose to express the reasoning in terms of games instead of
using trees, and that we needed the moves of Player 2 to include
the choice of
finite dimensional subspaces
$F_n$'s in which to pick the vectors $y_n$'s. This is due to
the fact that the basic sequence which witnesses the embedding of $X$
 into a given subspace
generated by a subsequence is not necessarily successive on the basis of $X$.

\

Let  $L,M$ be block-subspaces in ${\cal G}_{\Q}(X)$.
  Let $a \in
bb_{\Q}^{<\omega}(X)$ and $b \in (Fin_{\Q}(X)\times
{\Q}(X))^{<\omega}$ be such that
$|a|=|b|$ or $|a|=|b|+1$ (here $|x|$ denotes as usual the length of
the finite sequence $x$).
  Such a
couple
$(a,b)$ will be called a {\em state} of the
  game $G_{L,M}$ and the set of states will be written $St(X)$.
It is important to note that $St(X)$ is countable.
The empty sequence in $bb_{\Q}^{<\omega}(X)$ (resp.
$(Fin_{\Q}(X)\times {\Q}(X))^{<\omega}$) will be denoted by $\emptyset$.

\

  We define $G_{L,M}(a,b)$ intuitively as ``the game $G_{L,M}$
starting from the state $(a,b)$''.
More precisely, if $|a|=|b|$, then write $a=(a_1,\ldots,a_p)$ and
$b=(b_1,\ldots,b_p)$, with
$b_i=(B_i,\beta_i)$ for $i \leq p$.

  A $k$-th move for Player 1 is
$(x_k,n_k) \in {\Q}(X) \times \N$,
with $x_k \in L$, $x_1>a_p$ if $k=1$ and $a \neq \emptyset$, and
$x_k>x_{k-1}$ if $k \geq 2$.
  A $k$-th
move for Player 2 is
$(F_k,y_k) \in Fin_{\Q}(X) \times {\Q}(X)$ with $n_k \leq F_k \subset M$
and $y_k \in B_1 \oplus \ldots \oplus B_p \oplus F_1 \oplus \ldots
\oplus F_k$ for all $k$.

\[\begin{array}{cccccccc}
                   & I\ :& x_1,n_1 &      &  x_2,n_2 &     & \ldots &    \\
           G_{L,M}(a,b),\ \ &   &           &      &            &
&        &    \\
           |a|=|b|\ \ &   &           &      &            &     &        &    \\
                   & II:\ &          & F_1,y_1 &       & F_2,y_2 &   &
\ldots \end{array}\]

Player 2 wins $G_{L,M}(a,b)$ if the sequence
  $(\beta_1,\ldots,\beta_p,y_1,y_2,\ldots)$ is $C$-equi\-va\-lent to the
sequence $(a_1,\ldots,a_p,x_1,x_2,\ldots)$.

\

Now if $|a|=|b|+1$, then write $a=(a_1,\ldots,a_{p+1})$ and
$b=(b_1,\ldots,b_p)$, with
$b_i=(B_i,\beta_i)$ for $i \leq p$.

A first move for Player 1 is $n_1 \in \N$. A first move for Player 2 is
$(F_1,y_1) \in Fin_{\Q}(X) \times {\Q}(X)$ with $n_1 \leq F_1 \subset M$
and $y_1 \in B_1 \oplus \ldots \oplus B_p \oplus F_1$.

For $k \geq 2$, a $k$-th move for Player 1 is
$(x_k,n_k) \in {\Q}(X) \times \N$,
with $x_k \in L$, $x_2>a_{p+1}$ if $k=2$, and $x_k>x_{k-1}$ if $k >
2$; a $k$-th
move for Player 2 is
$(F_k,y_k) \in Fin_{\Q}(X) \times {\Q}(X)$ with $n_k \leq F_k \subset M$
and $y_k \in B_1 \oplus \ldots \oplus B_p \oplus F_1 \oplus \ldots \oplus F_k$.

\[\begin{array}{cccccccc}
                   & I\ :& n_1 &      &  x_2,n_2 &     & \ldots &    \\
           G_{L,M}(a,b),\ \ &   &           &      &            &
&        &    \\
           |a|=|b|+1\ \ &   &           &      &            &     &
&    \\
                   & II:\ &          & F_1,y_1 &       & F_2,y_2 &   &
\ldots \end{array}\]

Player 2 wins $G_{L,M}(a,b)$ if
  the sequence $(\beta_1,\ldots,\beta_p,y_1,y_2,\ldots)$ is
$C$-equi\-va\-lent to
the sequence $(a_1,\ldots,a_p,a_{p+1},x_2,\ldots)$.

\

We shall use the following classical stabilization process, called
"zawada" in \cite{P},
see also the proof by B. Maurey of Gowers' dichotomy theorem \cite{M}.
We define the following order relation on ${\cal G}_{\Q}(X)$: for
$M,N \in {\cal G}_{\Q}(X)$,
with $M=[m_i]_{i \in \N}, (m_i)_{i \in \N} \in bb_{\Q}(X)$, write
$M \subset^* N$ if there exists $p \in \N$ such
that $m_i \in N$ for all $i \geq p$.

Let
$\tau$ be a mapping defined on
${\cal G}_{\Q}(X)$ with values in the set $2^{\Sigma}$ of
subsets of some countable set $\Sigma$. Assume the map $\tau$ is
monotonous with
respect to $\subset^*$ on ${\cal G}_{\Q}(X)$
and to inclusion on $2^{\Sigma}$. Then by \cite{P} Lemma 2.1,  there exists a
block-subspace $M \in {\cal G}_{\Q}(X)$ which is stabilizing for
$\tau$, i.e. $\tau(N)=\tau(M)$
for every
$N
\subset^* M$.

\

We now define a map $\tau: {\cal G}_{\Q}(X) \rightarrow 2^{St(X)}$ by
$(a,b) \in \tau(M)$ iff  there exists $L \subset^* M$ such that Player
2 has a winning
   strategy for the game $G_{L,M}(a,b)$.

\begin{lemm} Let $M'$ and $M$ be in ${\cal G}_{\Q}(X)$. If
$M'
\subset^* M$ then $\tau(M') \subset
\tau(M)$.
\end{lemm}

\pf Let $M' \subset^* M$, let $(a,b) \in \tau(M')$, and let $L \subset^*
M'$ be such that
Player 2 has a winning strategy in
$G_{L,M'}(a,b)$ . Let $m$ be an
integer such that for any $x \in {\Q}(X)$, $x \in M'$ and
$\min(supp(x)) \geq m$ implies $x \in M$.
We describe a winning strategy for Player 2 in the game
$G_{L,M}(a,b)$: assume Player 1's $p$-th
move was $(n_p,x_p)$ (or just $n_1$ if it was the first move and $|a|=|b|+1$),
without loss of
generality
$n_p \geq m$. Let $(F_p,y_p)$ be the move prescribed
by the winning strategy for Player 2 in $G_{L,M'}(a,b)$. Then $F_p
\geq n_p \geq m$ and $F_p \subset
M'$ so $F_p \subset M$. The other conditions are satisfied to ensure
that we have described the $p$-th
move of a winning strategy for Player 2 in the game $G_{L,M}(a,b)$.
It remains to note that
$L \subset^* M$ as well to conclude that $(a,b) \in \tau(M)$.
\pff

\

By the stabilization lemma, there exists a block-subspace $M_0 \in
{\cal G}_{\Q}(X)$ such that for any $M
\subset^* M_0$,
$\tau(M)=\tau(M_0)$.

\

For $L,M \in {\cal G}_{\Q}(X)$ we shall write $L=^*M$ if $L \subset^* M$
and $M \subset^* L$.

We now define a map $\rho: {\cal G}_{\Q}(X) \rightarrow 2^{St(X)}$
  by
$(a,b) \in \rho(M)$ iff  there exists $L =^* M$ such that Player 2
has a winning
   strategy for the game $G_{L,M_0}(a,b)$.

\begin{lemm} Let $M'$ and $M$ be in ${\cal G}_{\Q}(X)$. If
$M'
\subset^* M$ then $\rho(M') \supset
\rho(M)$.
\end{lemm}

\pf Let $M' \subset^* M$, let $(a,b) \in \rho(M)$, and let $L =^* M$ be such that
Player 2 has a winning strategy in
$G_{L,M_0}(a,b)$.
Define $L'=M' \cap L$. As $L' \subset L$,
it follows immediately that Player 2 has a winning strategy in the game
$G_{L',M_0}(a,b)$.
It is also clear that $L'=^* M'$ so
$(a,b) \in \rho(M')$.
\pff

So there exists a block-subspace $M_{00} \in {\cal G}_{\Q}(X)$ of
$M_0$ which is stabilizing for
$\rho$, i.e. for any
$M
\subset^* M_{00}$,
$\rho(M)=\rho(M_{00})$.

\begin{lemm} $\rho(M_{00})=\tau(M_{00})=\tau(M_0).$
\end{lemm}
\pf First it is obvious by definition of $M_0$ that $\tau(M_{00})=\tau(M_0)$.

Let $(a,b) \in \rho(M_{00})$. There exists $L=^* M_{00}$ such that
Player 2 has a winning strategy in $G_{L,M_0}(a,b)$; as $L \subset^*
M_0$, this implies that
$(a,b) \in \tau(M_0)$.

Let $(a,b) \in \tau(M_{00})$. There exists $L \subset^* M_{00}$ such that
Player 2 has a winning strategy in $G_{L,M_{00}}(a,b)$. As $M_{00}
\subset M_0$, this is a winning
strategy for $G_{L,M_0}(a,b)$ as well. This implies that $(a,b) \in
\rho(L)$ and by the stablization
property for $\rho$, $(a,b) \in \rho(M_{00})$.
\pff

\

We now turn to the concluding part of the proof of Theorem
\ref{anyaisomorphic}.
By our assumption about $X$ just before Definition \ref{glm}, there
exists a block-sequence $(l_n)_{n
\in
\N}$ of
$bb_{\Q}(X)$ which is contained in $M_{00}$, and
$C$-continuously embeds into its subsequences, and without loss of
generality assume that $L_0:=[l_n, n \in \N]=M_{00}$. We fix an arbitrary
block-subspace $M$ of $L_0$ generated by a block-sequence in
$bb_{\Q}(X)$ and we shall prove that $L_0$ embeds into $M$.
By standard perturbation arguments this implies that $L_0$ is minimal.

  We construct by induction a subsequence $(a_n)_{n \in \N}$ of $(l_n)_{n
    \in \N}$, a sequence
$b_n=(F_n,y_n) \in (Fin_{\Q}(X)\times {\Q}(X))^{\omega}$ such that
$F_n \subset M$
and $y_n \in F_1 \oplus \ldots F_n$ for all $n \in \N$, and such that
   $((a_n)_{n \leq p},(F_n,y_n)_{n \leq p})
\in \rho(L_0)$ for all $p \in \N$.

  By Lemma \ref{firststepoftheinduction}, Player 2 has a winning
strategy in $G_{L_0,L_0}$, and so in particular $(\emptyset,\emptyset) \in
\rho(L_0)$ (recall that $\emptyset$ denotes the empty sequence in the
sets corresponding to the first
and second coordinates). This takes care of the first step of the induction.

Assume $(a,b)=((a_n)_{n \leq p-1},(F_n,y_n)_{n \leq p-1})$ is a state such that
$(a_n)_{n \leq p-1}$ is a finite subsequence of
$(l_n)_{n \in
\N}$, such that $F_n
\subset M$ and
$y_n
\in F_1
\oplus
\ldots F_n$ for all
$n \leq p-1$,
and such that $(a,b)
\in \rho(L_0)$.

  As $(a,b)$ belongs to
$\rho(L_0)$, there exists $L=^{*} L_0$ such that Player 2 has a winning
strategy in the game $G_{L,M_0}(a,b)$. In particular $L_0 \subset^* L$ so we may
 choose
$m_p$ large enough so that
$l_{m_p} > a_{p-1}$ and $l_{m_p} \in L$; we let Player 1 play $a_{p}=l_{m_p}$.
Player 2 has a winning strategy in the game
$G_{L,M_0}(a',b)$, where $a'=(a_n)_{n \leq p}$. In other words, $(a',b)$
belongs to $\rho(L_0)$. Now $\rho(L_0)=\tau(M)$, so there exists $L
\subset^* M$ such that Player 2
has a winning strategy in the game $G_{L,M}(a',b)$. Let Player 1 play
any integer
$n_p$, and $(F_p,y_p)$ with $F_p \subset M$ and $y_p \in F_1 \oplus \ldots
\oplus F_p$ be a move for Player 2 prescribed by that winning
strategy in response to $n_p$. Once again, Player 2 has a winning strategy
in $G_{L,M}(a',b')$, with $b'=(F_n,y_n)_{n \leq p}$; i.e. $(a',b') \in
\tau(M)=\rho(L_0)$.

\

To conclude, note that $(a_n,b_n)_{n \leq p} \in \rho(L_0)$ implies
  in particular that $(a_n)_{n \leq
p}$ and
$(y_n)_{n
   \leq p}$ are $C$-equivalent, and this is true for any $p \in \N$,
so $(a_n)_{n \in \N}$ and
$(y_n)_{n \in \N}$ are
$C$-equivalent. So $[a_n]_{n \in \N}$ $C$-embeds into
$M$.
Now $(a_n)_{n \in \N}$ is a subsequence of $(l_n)_{n \in \N}$ so by
our hypothesis, $L_0$ $C$-embeds
into
$[a_n]_{n \in \N}$ and thus $C^2$-embeds in $M$, and this concludes
the proof of Theorem
\ref{anyaisomorphic}.

\

 D. Kutzarova drew our attention to the dual $T^*$ of Tsirelson's space; it is minimal
\cite{CS} , but contains no block-minimal block-subspace (use e.g. \cite{CS} Proposition 2.4
and Corollary 7.b.3 in their $T^*$ versions, with Remark 1 after \cite{CS} Proposition
1.16).
 So Theorem \ref{anyaisomorphic} applies to
situations where the Theorem of Pelczar does not. On the other hand,
we do have:

\begin{coro} \label{coroblock} A Banach space with a Schauder basis which is saturated with
isomorphically homogeneous basic
sequences
   contains a block-minimal block-subspace.
\end{coro}

\pf Let $X$ have a Schauder basis and be saturated with isomorphically homogeneous basic
sequences. By the beginning of the
proof of Lemma \ref{seconddiagonalization}, we may assume $X$ is reflexive.
By Theorem \ref{anyaisomorphic}, there exists a minimal subspace $Y$
in
$X$, which is a block-subspace if you wish; passing
to a further
block-subspace assume
furthermore that $Y$ has an isomorphically homogeneous basis.
Take any block-subspace $Z$ of $Y=[y_n]_{n \in \N}$, then $Y$ embeds
into $Z$. By reflexivity and
standard perturbation results, some subsequence of $(y_n)_{n \in \N}$
spans a subspace which embeds
as a block-subspace of
$Z$. As
$(y_n)_{n \in \N}$ is isomorphically homogeneous,
this means that $Y$ embeds as a block-subspace of $Z$. \pff

\

We recall that a Banach space is said to be ergodic if the relation
$E_0$ is Borel reducible to the
relation of isomorphism between its subspaces.

\begin{coro} \label{ergodicorminimal} A Banach space is either
ergodic or contains a minimal subspace.
\end{coro}

\pf We prove the stronger result that if $X$ is a Banach space with a
Schauder basis, then either
$E_0$ is Borel reducible to isomorphism between block-subspaces of
$X$ or $X$ contains a
block-minimal block-subspace.

Assume $E_0$ is not Borel reducible to isomorphism between
block-subspaces of $X$. By
\cite{R}, Theorem 19,  any block-sequence in $X$ has an
isomorphically homogeneous subsequence. In particular $X$ is saturated with
isomorphically homogeneous sequences, so apply Corollary
\ref{coroblock}.
\pff

\begin{coro} A Banach space $X$ contains a minimal subspace or the
relation $E_0$ is Borel reducible
to the relation of biembeddability between subspaces of $X$.
\end{coro}

\pf Note that the relation $\sim^{emb}$ of biembeddability between
subspaces of $X$
is analytic. By
\cite{R} Theorem 15, if
$E_0$ is not Borel reducible to biembeddability between subspaces of
$X$, then every basic sequence in $X$ has a subsequence $(x_n)_{n \in
\N}$ which is homogeneous for
the relation between subsequences corresponding to $\sim^{emb}$, that is,
for any subsequence $(x_n)_{n \in I}$ of $(x_n)_{n \in \N}$,
$[x_n]_{n \in I} \sim^{emb} [x_n]_{n
\in \N}$. This means that $(x_n)_{n \in \N}$ embeds into its subsequences.
  So $X$ is saturated with basic sequences
which embed into their subsequences.\pff

\

We conclude this section with a remark about the proof of Theorem \ref{anyaisomorphic}.
The
sequences
$(m_p)_{p
\in
\N}
\in [\omega]^{\omega}$ (associated to a subsequence of
$(l_n)_{n \in \N}$) and $b_p=(F_p,y_p) \in (Fin_{\Q}(X)\times {\Q}(X))^{\omega}$ (with
$(y_p)_{p \in \N}$ $C$-equivalent to $(l_{m_p})_{p \in \N}$) in our final induction may
clearly be chosen with
$F_p
\subset M_p$ for all
$p$, for an arbitrary sequence
$(M_p)_{p \in \N}$ of block-subspaces of $L_0$.
Also, $(l_n)_{n \in \N}$ $C$-continuously embeds into its subsequences, i.e. there is a
continuous map $f:[\omega]^{\omega} \rightarrow bb_{\Q}(X)$ such that $f(A)$ is
$C$-equivalent to $(l_n)_{n \in \N}$ for all $A \in bb_{\Q}(X)$.

By combining these two facts, it is easy to see that Player 2 has a winning strategy
to produce a sequence $(y_n)_{n \in \N}$ which is $C^2$-equivalent
to $(l_n)_{n \in \N}$, in a
"modified" Gowers' game, where a $p$-th move for Player 1 is a
block-subspace $Y_p \in {\cal G}_{\Q}(X)$, with $Y_p \subset L_0$, and a $p$-th move for
Player 2 is a couple $(F_p,y_p) \in (Fin_{\Q}(X)\times {\Q}(X))^{\omega}$ with $F_p \subset
Y_p$ and
$y_p
\in F_1 \oplus \ldots \oplus F_p$. 

This is an instance of a result with a Gowers' type game
where Player 2 is allowed to play sequences of vectors which are not necessarily
block-basic sequences. 

\section{Topological 0-1 laws for block-sequences.}

In this section $X$ denotes a Banach space with a normalized basis
$(e_n)_{n \in \N}$.
It will be
necessary to restrict our attention to normalized block-bases in $X$
to use compactness properties.
We denote by $bb(X)$ the set of normalized block-bases on $X$.
Let  $Q(X)$ be the set
of normalized blocks of the basis that are a multiple of some block with
rational coordinates (or coordinates in $\Q+i\Q$ in the complex case).
  We
denote by $bb_d(X)$ the set of block-bases of vectors in $Q(X)$ (here
"$d$" stands for "discrete",
this notation was introduced in \cite{FR1}).
We consider $bb_d(X)$ as a topological space,
equipped with the product topology
of the discrete topology on $Q(X)$, which turns it
  into a Polish space.

The notation $bb_d^{<\omega}(X)$  will denote
the set of finite block-sequences with blocks in $Q(X)$.
For a finite block
sequence
$\tilde{x}=(x_1,\ldots,x_n) \in bb_d^{<\omega}(X)$, we denote by
$N(\tilde{x})$
the set of elements of $bb_d(X)$ whose first $n$ vectors are
$(x_1,\ldots,x_n)$; this is the basic open set associated to $\tilde{x}$.

If $s$ is a finite block-basis and $y$ is a finite or infinite block-basis
supported after $s$, denote by $s^{\frown}y$  the concatenation of
$s$ and $y$. The notation $x=(x_n)_{n \in \N}$ will be reserved to denote
an infinite block-sequence, and $[x]$ will denote its closed linear
span; $\tilde{x}$ will denote a
finite block-sequence, and $|\tilde{x}|$ its
length as a sequence, $supp(\tilde{x})$ the union of the supports of the
terms of $\tilde{x}$. For two finite block-sequences
$\tilde{x}=(x_1,\ldots,x_n)$
and
$\tilde{y}=(y_1,\ldots,y_m)$, write $\tilde{x}<\tilde{y}$ to mean that they are
successive, i.e. $x_n<y_1$. For a sequence of successive finite block-sequences
$(\tilde{x}_i)_{i
\in I}$, we denote the concatenation of the block-sequences by
$\tilde{x}_1^{\frown} \ldots^{\frown} \tilde{x}_n$ if the sequence is finite with
$I=\{1,\ldots,n\}$, or
$\tilde{x}_1^{\frown} \tilde{x}_2^{\frown}\ldots$ if it is infinite, and
we denote by $supp(\tilde{x}_i, i \in I)$ the support of the
concatenation, by
$[\tilde{x}_i]_{i
\in I}$ the closed linear span of the concatenation.

\

\subsection{A principle of topological 0-1 law for
block-sequences.}

We recall a characterization of comeager  subsets of $bb_d(X)$ which was
proved in \cite{FR2}.
If $A$ is a subset of $bb_d(X)$ and $\Delta=(\delta_n)_{n \in \N}$ is a
sequence of
  positive real numbers, we denote by $A_{\Delta}$ the
{\em $\Delta$-expansion of
$A$} in $bb_d(X)$, that is $x=(x_n) \in A_{\Delta}$ iff
there exists $y=(y_n) \in A$ such that $\norm{y_n-x_n} \leq
\delta_n, \forall n \in
\N$. Such an $y$ will be called a {\em $\Delta$-perturbation of $x$}.
Given a finite block-sequence $\tilde{x}=(x_1,\ldots,x_n)$, we say that a
(finite or infinite)
block-sequence $(y_i)_{i \in \N}$ {\em passes through} $\tilde{x}$ if
there exists some integer $m$ such that $\forall 1 \leq i \leq n$,
$y_{m+i}=x_i$.

\begin{prop}(V. Ferenczi, C. Rosendal \cite{FR2}) \label{comeager}
Let $X$ be a Banach space with a Schauder basis. Let $A$ be comeager  in $bb_d(X)$. Then
  for all $\Delta>0$, there exists a sequence $(\tilde{a}_n)_{n \in \N} \in
(bb_d^{<\omega}(X))^{\omega}$ of suc\-ces\-sive fini\-te
block-se\-quences
   such that any block-sequence of $bb_d(X)$ passing trough
   infinitely many of the $\tilde{a}_n$'s is in $A_{\Delta}$.
   \end{prop}

As was noted in \cite{FR2}, the property in the conclusion of this
proposition  is essentially (i.e. up to perturbation) a
characterization of comeager sets in $bb_d(X)$. Indeed, it easily implies
that $A_{\Delta}$ is comeager.

Let $A$ have the Baire Property, that is, there exists an open set $U$ such
that
$A \setminus U$ and $U \setminus A$ are meager. Then
either $A$
is meager, or $A$ is comeager in $N(\tilde{x}_0)$ for some finite
block-sequence $\tilde{x}_0 \in bb_d^{<\omega}(X)$.
This fact is to be combined with Proposition
\ref{comeager}. For example, $A$ is comeager in $N(\tilde{x}_0)$ will
imply that for all $\Delta>0$, there exists a sequence of successive
finite block-sequences
   $(\tilde{a}_n)_{n \in \N}$ such that any element of $bb_d(X)$ passing trough
$\tilde{x}_0$ and  infinitely many of the $\tilde{a}_n$'s is in $A_{\Delta}$.
This result will take more interest if one assumes a few natural additional
properties for the set $A$.

\

In the following we identify a block-sequence $(x_k)_{k \in K}$
indexed on some infinite subset
$K=\{k_1,k_2,\ldots\}$ of
$\N$ (where $(k_n)_{n \in \N}$ is increasing), with the block-sequence
$(x_{k_n})_{n \in \N}$ indexed on $\N$; so given an infinite block-sequence,
we may always chose the more convenient way to index it. We do the similar
identification for finite block-sequences.

For $(x_n)_{n \in \N}$ a block-sequence and $n_0 \in \N$, we call
{\em $n_0$-modification of
$(x_n)_{n \in \N}$} a block-sequence $(y_n)_{n \in \N}$ such that
$x_n=y_n$ for all $n > n_0$. An $n_0$-modification of $(x_n)_{n \in \N}$ for
some
$n_0$ will be
called a {\em finite modification} of $(x_n)_{n \in \N}$.
For a block-sequence $(x_n)_{n \in \N}$ of
$X$, a couple $((x_n)_{n \in I},(x_n)_{n \in J})$ of block-sequences
associated to a partition of $\N$ in two infinite sets $I$ and $J$ will be
called a {\em partition} of $(x_n)_{n \in \N}$.

Related to the notion of support is the useful notion of
range: the {\em range} $ran(x_0)$ of $x_0 \in X$ is the smallest interval of
integers containing the support of $x_0$. If $x=(x_n)_{n \in I}$ is a finite or
infinite  block-sequence, $ran(x)$ will denote the union
$\cup_{n \in I} ran(x_n)$. When $x=(x_n)_{n \in I}, y=(y_n)_{n \in
J}$ are finite or infinite block-sequences
whose ranges are disjoint, we call {\em concatenation of $x$ and $y$} the
unique (up to the choice of $K$) block-sequence $z=(z_n)_{n \in K}$ such that
  $\{z_n, n \in K\}=\{x_n, n \in I\} \cup \{y_n, n \in J\}$.

\

We are now ready to state our principle of
topological 0-1 law for block-sequences.

\begin{theo}\label{topologicallaw} (Topological 0-1 law for block-sequences)
Let $X$ be a Banach space with a Schauder basis. Assume $A \subset bb_d(X)$
has the Baire Property and is invariant
by finite modifications. Then $A$ is either meager or comeager in $bb_d(X)$.
\end{theo}

This is a corollary of the following quantified version:

\begin{prop}\label{quantifiedtopologicallaw} Let $X$ be a Banach
space with a Schauder basis.
Let $(A_N)_{N \in \N}$ be an increasing sequence of  subsets of
$bb_d(X)$ with the Baire Property, and let
$A=\cup_{N \in \N} A_N$. Assume
that for any $N \in \N$
and $n_0 \in \N$, there exists $K(N,n_0) \in \N$ such that
whenever $(x_n)_{n \in \N}$ belongs to $A_N$, then any $n_0$-modification
of $(x_n)_{n \in \N}$ belongs to $A_{K(N,n_0)}$.

Then either $A$ is meager in $bb_d(X)$,
either there exists $K \in \N$ such that $A_K$ is comeager in $bb_d(X)$.
\end{prop}

\pf
We assume $A$ is non-meager, then for some $N \in \N$, $A_N$ is
non-meager. We reproduce a proof of \cite{FR2}. As $A_N$ has the Baire
property, it is comeager in some basic open set $U$, of the form
$N(\tilde{x})$, for some finite block-sequence $\tilde{x} \in
bb_d^{<\omega}(X)$.

We now prove that $A_K$ is comeager in $bb_d(X)$ for
$K=Nc(2\max(supp(\tilde{x})))$ (for $n \in \N$, $c(n)$ denotes
an integer such for any Banach space $X$, any $n$-codimen\-sional
subspaces of $X$ are $c(n)$-isomorphic, see e.g.\cite{FR1}). So
let us assume $V=N(\tilde{y})$ is some basic open set in $bb_d(X)$
such that $A_K$ is meager in $V$. We may assume that
$|\tilde{y}|>|\tilde{x}|$ and write
$\tilde{y}=\tilde{x}^{\prime\frown}\tilde{z}$ with
$\tilde{x}<\tilde{z}$ and $|\tilde{x}^{\prime}| \leq
\max(supp(\tilde{x}))$. Choose $\tilde{u}$ and $\tilde{v}$ in
$bb_d^{<\omega}(X)$ such that
$\tilde{u},\tilde{v}>\tilde{z}$, $|\tilde{u}|=|\tilde{x}^{\prime}|$ and
$|\tilde{v}|=|\tilde{x}|$, and such that
$\max(supp(\tilde{u}))=\max(supp(\tilde{v}))$. Let $U'$ be the
basic open set $N(\tilde{x}^{\frown}\tilde{z}^{\frown}\tilde{u})$
and let $V'$ be the basic open set
$N(\tilde{x}^{\prime\frown}\tilde{z}^{\frown}\tilde{v})$. Again $A_N$
is comeager in $U'$ while $A_K$
is meager in $V'$.

Now let $T$ be the canonical map from $U'$ to $V'$.
For all $u$ in $U'$, $T(u)$ differs from at most
  $|\tilde{x}|+\max(supp(\tilde{x})) \leq
2\max(supp(\tilde{x}))$ vectors from $u$, so $[T(u)]$ is
$c(2\max(supp(\tilde{x})))$ isomorphic to $[u]$. Since
$K=Nc(2\max(supp(\tilde{x})))$ it follows that $A_K$ is
comeager  in $V' \subset V$. The contradiction follows by choice
of $V$.
\pff

Proposition \ref{comeager}
characterizes  meager and comeager sets in the conclusion of Theorem
\ref{topologicallaw}.
This leads us to the following theorem.

\begin{theo}\label{dichotomytheorem}
Let $X$ be a Banach space with a Schauder basis.
Let $A$ be a subset of $bb_d(X)$ with the Baire Property, which
  is stable by $\Delta$-perturbations for some
$\Delta>0$, by finite modifications, and by taking subsequences.
Assume that any sequence $(\tilde{x}_n)_{n \in \N} \in
(bb_d^{<\omega}(X))^{\omega}$ of successive
finite block-sequences
   admits a subsequence $(\tilde{x}_{n_k})_{k \in \N}$ such that the
block-sequence
$\tilde{x}_{n_1}^{\frown}\tilde{x}_{n_2}^{\frown}\ldots$ belongs to $A$.

Then every block-sequence in $bb_d(X)$ admits a partition in a couple
of elements
of $A$.

If furthermore, the set $A$ is stable by  conca\-te\-na\-tion of
pairs of block-se\-quen\-ces,
then
$bb_d(X)=A$.

\end{theo}

Once again this is a corollary of  a quantified version:

\begin{prop}\label{theo3}
Let $X$ be a Banach space with a Schauder basis.
Let $(A_N)_{N \in \N}$ be an increasing sequence of subsets of
$bb_d(X)$ with the Baire Property,
such that:

\

(a) there exists $\Delta>0$ such that for any $N \in \N$,
  there exists $K_1(N) \in \N$ such that $(A_N)_{\Delta} \subset A_{K_1(N)}$.

(b) for any $N \in \N$
and $n_0 \in \N$, there exists $K_2(N,n_0) \in \N$ such that
whenever $(x_n)_{n \in \N}$ belongs to $A_N$, then any $n_0$-modification
of $(x_n)_{n \in \N}$ belongs to $A_{K_2(N,n_0)}$.

(c) for any $N \in \N$, there exists $K_3(N) \in \N$ such that
whenever $(x_n)_{n \in \N}$ belongs
  to $A_N$ then any subsequence of $(x_n)_{n \in \N}$  belongs to $A_{K_3(N)}$.

\

Let $A=\cup_{N \in \N}A_N$.
Assume that  any sequence $(\tilde{x}_n)_{n \in \N} \in
(bb_d^{<\omega}(X))^{\omega}$ of successive
finite block-sequences
admits a subsequence $(\tilde{x}_{n_k})_{k \in \N}$ such that
the block-sequence
$\tilde{x}_{n_1}^{\frown}\tilde{x}_{n_2}^{\frown}\ldots$ belongs to $A$.

Then there exists $N \in \N$ such that every block-sequence
  in $bb_d(X)$ has a partition  in two elements
  of $A_N$.
If furthermore,

\

(d) for any $N \in \N$, there exists $K_4(N) \in \N$ such that any
concatenation
of a couple of block-sequences in $A_N^2$ belongs to $A_{K_4(N)}$,

\

then $bb_d(X)=A_N$ for some $N \in \N$.
\end{prop}

\pf The part which is a consequence of (d) is obvious once we prove the first
part of the proposition.
We note that by Proposition \ref{quantifiedtopologicallaw}, or $A$ is meager,
or $A_N$ is comeager for some $N \in \N$.
By (a), there is some $\Delta>0$ such that $A=A_{\Delta}$.  It follows that
$A_{\Delta} \cap A^C=\emptyset$, that is $(A^C)_{\Delta} \cap A = \emptyset$.

If $A$ is meager, Proposition \ref{comeager} gives us a sequence of successive
finite block-sequences $(\tilde{x}_n)_{n \in \N}$ such that, in particular,
$\tilde{x}_{n_1}^{\frown}\tilde{x}_{n_2}^{\frown}\ldots$ is in
$(A^C)_{\Delta}$ for every subsequence $(\tilde{x}_{n_k})_{k \in \N}$,
So for no subsequence $(\tilde{x}_{n_k})_{k \in \N}$,
$\tilde{x}_{n_1}^{\frown}\tilde{x}_{n_2}^{\frown}\ldots$ is in $A$.

So $A_N$ is comeager for some $N \in \N$. Applying Proposition
\ref{comeager}, and up
to modifying $N$, let
   $(\tilde{a}_n)_{n \in \N}$ be a sequence of successive block-sequences such
   that every block-sequence passing through infinitely many of the
   $\tilde{a}_n$'s is in $A_N$.

Let now $(x_n)_{n \in \N}$ be an arbitrary block-sequence in $bb_d(X)$.
We note that we may find a partition of $(x_n)_{n \in \N}$ in two subsequences
$(x_n)_{n \in I}$ and $(x_n)_{n \in J}$, and a subsequence
$(\tilde{a}_{n_k})_{k \in \N}$
of
$(\tilde{a}_n)_{n \in \N}$ such that
  $(x_n)_{n \in I}$ and $(\tilde{a}_{n_{2k}})_{k \in \N}$ have disjoint
ranges (let $(i_n)_{n \in \N}$ denote their concatenation)
and such that $(x_n)_{n \in J}$ and $(\tilde{a}_{n_{2k-1}})_{k \in
\N}$ have disjoint
ranges (let $(j_n)_{n \in \N}$ denote their concatenation).

Now $(i_n)_{n \in \N}$ belongs to $A_N$, so by (c), for some $N' \in
\N$, $(x_n)_{n \in I}$
belongs to $A_{N'}$, and likewise $(x_n)_{n \in J}$ belongs to
$A_{N'}$.
\pff

\

In particular we deduce a uniformity principle from Proposition
\ref{theo3}. Under its
hypotheses, and if every block-sequence of $bb_d(X)$
  is in $A$, then there exists $N \in \N$ such that every block-sequence of
$bb_d(X)$ is in $A_N$.
This method was first used in \cite{FR2} to study the property of
complementable embeddability
(\cite{FR2} Proposition 17).

Before passing to applications,
we note that there is a case when Proposition \ref{theo3} is not so
interesting. It is when the set
$A$ is $F_{\sigma}$ (in which case we may assume $A$ is the union of an
increasing sequence of closed sets $(A_n)_{n \in \N}$). In that case,
there is a much more direct
proof of it, which does not use the subsequence hypothesis (c) nor the
concatenation hypothesis (d).
A typical instance of this situation is when $A$ is the set of
block-sequences in
$bb_d(X)$ which are equivalent to a given basic sequence.

\begin{rema} Let $X$ be a Banach space with a Schauder basis.
Let $(A_n)_{n \in \N}$ be an increasing sequence of closed subsets of
$bb_d(X)$ and let
$A=\cup_{n \in \N} A_n$. Assume hypotheses (a) and (b) from Proposition
\ref{theo3} are satisfied.
  If
  any sequence $(\tilde{x}_n)_{n \in \N} \in
(bb_d^{<\omega}(X))^{\omega}$ of successive finite
block-sequences
   admits a subsequence $(\tilde{x}_{n_k})_{k \in \N}$ such that the
block-sequence
$\tilde{x}_{n_1}^{\frown}\tilde{x}_{n_2}^{\frown}\ldots$ belongs to $A$, then
$bb_d(X)=A_N$ for some $N \in \N$.
\end{rema}

\pf Assume the conclusion does not hold.
Note that for every $k \in \N$, every finite block-sequence
$\tilde{x}_0 \in bb_d^{<\omega}(X)$,
  there
exists
a finite block-sequence $\tilde{x} \in bb_d^{<\omega}(X)$ such that
$\tilde{x}_0^{\frown} \tilde{x}$
is not extendable in an element of $A_k$. Otherwise, by closedness of
$A_k$, we would
have that $\tilde{x}_0^{\frown} x \in A_k$ for all $x \in bb_d(X)$
supported after
$\tilde{x}_0$, and using (b) we would deduce that every $x \in bb_d(X)$
belongs to $A_K$ for some $K$ depending on $k$ and the length of the support
of
$\tilde{x}_0$.

It follows easily that for any $k \in \N$, for any finite block-sequences
$\tilde{x}_0, \ldots, \tilde{x}_n$ in $bb_d^{<\omega}(X)$, there
exists $\tilde{x}$ such that
$\tilde{x}_i^{\frown} \tilde{x}$ is extendable in an element of $A_k$
for no $i \leq n$.

Using this fact, we find $\tilde{a}_1 \in bb_d^{<\omega}(X)$ not
extendable in $A_1$, and by
induction, for any $n \in \N$, a finite block-sequence
$\tilde{a}_n>\tilde{a}_{n-1}$ such that
$\tilde{a}_{i_1}^{\frown}\ldots^{\frown}\tilde{a}_{i_p}^{\frown}\tilde{a}_n$
is extendable in $A_n$ for no finite sequence $i_1<\cdots<i_p<n$.

Consider now the sequence $(\tilde{a}_n)_{n \in \N}$. By construction,
for any subsequence $(\tilde{a}_{n_k})_{k \in \N}$, the block-sequence
$\tilde{a}_{n_1}^{\frown}\tilde{a}_{n_2}^{\frown}\ldots$ is in $A_{n_k}$ for
no $n_k, k \in \N$, so does not belong to $A$, a contradiction with
the hypotesis.
\pff

\subsection{Isomorphism between block-subspaces.}

Recall that two basic sequences $(x_n)_{n \in \N}$ and $(y_n)_{n \in \N}$ are said to be
{\em permutatively equivalent} if there is a permutation $\sigma$ on $\N$ such that
$(x_n)_{n
\in
\N}$ is equivalent to $(y_{\sigma(n)})_{n \in \N}$. By Zippin's theorem and by a result of
Bourgain, Casazza, Lindenstrauss and Tzafriri
\cite{BCLT}, the homogeneity problem is solved for equivalence or permutative equivalence
of normalized block-sequences, the only solutions being the canonical bases
of $c_0$ or
$l_p$, $1 \leq p<+\infty$. Unfortunately, the situation is not nearly as nice when we
replace permutative equivalence by isomorphism. We shall say that a Banach space $X$ with a
basis is {\em block-homogeneous} if
$X$ is isomorphic to all its block-subspaces. We recall our conjecture:

\begin{conj}\label{blockhomogeneous} Let $X$ be a Banach space with an
   (uncon\-ditional) Schau\-der basis
$(e_n)_{n \in \N}$, which is block-homo\-geneous. Does it
follow that
$X$ is isomorphic to $c_0$ or $l_p$?
\end{conj}

Note that an (even unconditional) basis $(e_n)_{n \in \N}$ of a
block-homogeneous Banach
space $X$  need not be equivalent to the canonical basis of $c_0$ or
$l_p$: for $1<p<+\infty$, $X=(\oplus_{n \in \N} l_2^n)_p$
(with the associated canonical
basis) is isomorphic to $l_p$, and every block-subspace of $X$
is complemented in $X$ (\cite{LT}, Proposition 2.a.12) and thus
isomorphic to $l_p$ as well.

\

In some special cases however, results of uniqueness of unconditional basis
will allow us to
pass from isomorphism to permutative equivalence and use the previous results.
We recall some definitions and results from \cite{CK2}. A sequence
space $X$ is said to be
{\em left (resp. right) dominant} if there exists a constant $C \geq
1$ such that
whenever $(u_i)_{i \leq n}$ and $(v_i)_{i \leq n}$ are finite block-sequences,
with $\norm{u_i} \geq \norm{v_i}$ (resp. $\norm{u_i} \leq \norm{v_i}$) and
$v_i>u_i$ for all $i \leq n$, then
$\norm{\sum_{i=1}^n v_i} \leq C\norm{\sum_{i=1}^n u_i}$
(resp. $\norm{\sum_{i=1}^n u_i} \leq C\norm{\sum_{i=1}^n v_i}$.
When $X$ is left or right dominant, then there exists exactly one $r=r(X)$ such
that
$l_r$ is finitely disjointly representable in $X$, and we call $r$ the index
of $X$.

We refer to \cite{LT}, \cite{K} for the definition of and background
about Banach lattices.
If $X$ and $Y$ are Banach lattices, a bounded linear operator $V:X \rightarrow
Y$ is called a {\em lattice homomorphism} if $V(x_1 \vee x_2)=Vx_1 \vee Vx_2$
for all $x_1,x_2 \in X$.
  Following \cite{CK2}, define a Banach lattice $X$ to
be {\em sufficiently lattice-euclidean} if there exists $C \geq 1$ such thar
for all $n \in \N$, there exist operators $S:X \rightarrow l_2^n$ and $T:l_2^n
\rightarrow X$ such that $ST=I_{l_2^n}$, $\norm{S}\norm{T} \leq C$ and such
that
$S$ is a lattice homomorphism.
This is equivalent to saying that $l_2$ is finitely representable as a
complemented sublattice of $X$. A Banach lattice which is not sufficiently
lattice-euclidean is said to be {\em anti-lattice euclidean}.

For an unconditional basis $(x_n)_{n \in \N}$ of a Banach space (seen as a
Banach lattice), being
sufficiently lattice-euclidean is the same as
  having, for some $C \geq 1$ and every $n \in \N$,
a $C$-complemented, $C$-isomorphic copy of $l_2^n$
  whose basis is disjointly supported on $(x_n)_{n \in \N}$.

\begin{prop}\label{Schlum} Let $X$ be a Banach space with a normalized unconditional basis
$(e_n)_{n \in \N}$ which is isomorphically homogeneous.
Assume $(e_n)_{n \in \N}$ is
   right or left dominant with $r(X) \neq 2$ and that $(e_n)_{n \in \N}$
   is equivalent to $(e_{2n})_{n \in \N}$. Then
$(e_n)_{n \in \N}$ is equivalent to the
canonical
   basis of $l_p, p \neq 2$ or $c_0$.
\end{prop}

\pf Let $(y_n)_{n \in \N}$ be any subsequence of $(e_n)_{n \in \N}$, and
 $Y=[y_n]_{n \in \N}$. The sequence $(y_n)_{n \in \N}$ is equivalent
 to an unconditional basis $(u_n)_{n \in \N}$
of
$X$.
 It is enough to note now that the proof of \cite{CK2}, Theorem 5.7, is still valid as long
as we prove
that $(u_n)_{n \in \N}$ is anti-lattice
euclidean. But this is clear because $r(Y)=r(X) \neq 2$. So
$(y_n)_{n \in \N}$ must be permutatively equivalent to $(e_n)_{n \in \N}$.

It follows by \cite{BCLT} Proposition 6.2 that some subsequence $(v_n)_{n \in \N}$ of
$(e_n)_{n \in \N}$ is subsymmetric. By \cite{CK2}, $X$ is asymptotically $c_0$ or $l_p$ for
some $p \neq 2$, so $(v_n)_{n \in \N}$ is equivalent to the canonical basis of $c_0$
or $l_p$, $p \neq 2$ and $(e_n)_{n \in \N}$ as well.
\pff

\

The right or left dominant hypothesis in Proposition \ref{Schlum} cannot be removed:
the canonical basis $(e_n)_{n \in \N}$ of Schlumprecht's space $S$ \cite{S}
 is unconditional,
subsymmetric, but $S$
does not even contain a copy of $c_0$ or $l_p$.

 It is of interest to note that $S$ is however quite homogeneous in some sense:
any  constant coefficient block-subspace of $S$ is
 isomorphic to $S$ (see \cite{KL}, Remark before Proposition 9, for the proof and \cite{LT}
for the definition). So $S$ is an example of a non
$c_0$ or
$l_p$, yet "constant coefficient block-homogeneous" sequence space. This contrasts with the
 Theorem of Zippin (resp. the Theorem of Bourgain, Casazza, Lindenstrauss,
Tzafriri) for equivalence (resp. permutative equivalence) which can be proved using only
 constant coefficient block-sequences in $X$ (\cite{BCLT}).

\

The question of uniformity in the homogeneous Banach space problem was raised
by Gowers \cite{G}. Of course, since a homogeneous Banach space must be
isomorphic to $l_2$, it is
trivial that if $X$ is homogeneous, then there exists a constant $C \geq 1$
such that $X$ is $C$-isomorphic to any of its subspaces. However, there does
not seem to be a direct proof of this fact. Note also that uniformity is the
first step in the proof of the theorem of Zippin. So the following
question is natural:

\begin{conj} \label{uniformlyblockhomogeneous} Let $X$ be a Banach
space with an (unconditional) basis
   $(e_n)_{n \in \N}$. Assume $X$ is block homogeneous. Does
   there exists $C \geq 1$ such that  $X$ is $C$-block homogeneous?
\end{conj}

By a $C$-block homogeneous Banach space with a basis,
  we mean a Banach space  $C$-isomorphic
to all its block-subspaces.

\

As a partial result, we may use the primeness of the spaces $c_0$ and
$l_p$ to get a
positive answer to Question
\ref{uniformlyblockhomogeneous} when $X$ is isomorphic to $l_p$ or $c_0$:

\begin{prop} \label{isomorphictolp} Let $p \geq 1$. Let $X$ be a Banach
  space with an unconditional basis. Assume
  any sequence
$(\tilde{x}_n)_{n \in \N} \in (bb_d^{<\omega}(X))^{\omega}$
of successive finite block-sequences
   admits a subsequence $(\tilde{x}_{n_k})_{k \in \N}$ such that the
block-subspace
$[\tilde{x}_{n_1}^{\frown}\tilde{x}_{n_2}^{\frown}\ldots]$ is isomorphic
  to $l_p$.
Then $X$ is isomorphic to $l_p$ , and furthermore,
there exists $C \geq 1$, such that all block-subspaces of $X$ are
   $C$-isomorphic to $l_p$.
  The similar result holds for $c_0$.
\end{prop}

\pf We may assume the unconditional basis of $X$ is $1$-unconditional (then all
canonical projections on subspaces spanned by subsequences are of norm $1$).
  The set $A_N=\{(x_n)_{n \in \N} \in bb_d(X): [x_n]_{n \in \N}
\simeq^N l_p\}$ is
analytic and so has Baire Property (this is true of any isomorphism
class in $bb_d(X)$, see
\cite{FR2} about this). We check the hypotheses of Proposition
\ref{theo3}. Given $\epsilon>0$, there
exists $\Delta>0$ such that
the $\Delta$-perturbation of a block-sequence $(x_n)_{n \in \N}$ in $bb_d(X)$
spans a space which is $1+\epsilon$ isomorphic to $[x_n]_{n \in \N}$, so (a)
follows.
(b) is true with $K(N,n_0)=Nc(n_0)$ (here $c(n)$ is the previously
used constant
such that in any Banach space, any two subspaces of codimension $n$ are
$c(n)$-isomorphic). If $[x_n]_{n \in \N}$ is $C$-isomorphic
to $l_p$, and if $(x_{n_k})_{k \in \N}$ is a subsequence of $(x_n)_{n \in \N}$,
then as $[x_{n_k}]$ is $1$-complemented in $[x_n]_{n \in \N}$, it is
$C$-isomorphic to a $C$-complemented subspace of $l_p$, so is $K(C)$-isomorphic
to $l_p$, for some constant $K(C)$.
Finally it is easy to check that if $x,y$ in $bb_d(X)$ are disjointly
supported, and $[x]$ and $[y]$
are
$C$-isomorphic to
$l_p$, then the concatenation of $x$ and $y$ will span a subspace
which is $k(C)$ isomorphic
to
$l_p$, for some constant $k(C)$.
\pff

\

In particular, if $X$ has a block-homogeneous unconditional basis and is
isomorphic to $l_p$ or $c_0$, then it is $C$-block-homogeneous for some $C
\geq 1$. 

 To conclude this section, it is worth
noting the form that our topological 0-1 law  takes when
$A$ is really an isomorphic
property of the span of a block-sequence in $bb_d(X)$.

\begin{theo}
Let $X$ be a Banach space with an unconditional basis $(e_n)_{n \in \N}$.
Let $P$ be an isomorphic property of Banach spaces such that
$A=\{(x_n)_{n \in \N} \in
bb_d(X): [x_n]_{n \in N}\ has\ P\}$ has Baire Property, and which is stable
  by taking complemented subspaces.

Assume  any sequence of successive finite block-sequences
$(\tilde{x}_n)_{n \in \N}$
   admits a subsequence $(\tilde{x}_{n_k})_{k \in \N}$ such that the
block-subspace
$[\tilde{x}_{n_1}^{\frown}\tilde{x}_{n_2}^{\frown}\ldots]$ satisfies $P$.

Then every block-subspace of $X$ is the sum of two disjointly supported
block-subspaces
satisfying $P$.

Assume furthermore that
any direct sum of two spaces with $P$ satisfies $P$,
then every block-subspace
in $bb_d(X)$ satisfies $P$.
\end{theo}

\section{Topological 0-1 law for subspaces with a successive
finite dimensional decomposition.}

We now turn to subspaces of a space $X$ with a Schauder basis, which have a
successive finite dimensional decomposition on the basis. There is a
natural discretization of the set of such spaces, introduced in \cite{FR2}.

   We
  say that two finite-dimensional subspaces $F$ and $G$ of $X$
are {\em successive}, and write $F<G$, if they are different from $\{0\}$ and
for any $0 \neq x \in
F$,
$0 \neq y
\in G$,
$x$ and $y$ are successive.
A space with a {\em successive finite dimensional decomposition (or successive
  FDD) in $X$} is a
subspace of $X$ of the form
$\oplus_{k \in \N} F_k$,
with successive, finite-dimensional subspaces $F_k$. The associated sequence
$(F_k)_{k \in \N}$ will be called a {\em sequence of successive finite
  dimensional subspaces}.  Such a sequence {\em passes
through} a finite sequence of successive finite dimensional subspaces
$(A_i)_{1 \leq i \leq I}$ if there exists $k$ such that
$F_{k+i}=A_i$ for all $1 \leq i \leq I$.
If the sequence $(A_i)_i$ is a length $1$ sequence $(A)$,
we shall just say that $(F_k)_{k \in \N}$ passes through
$A$.

We let $fdd(X)$ be the set of infinite sequences of successive
finite-dimensional subspaces, and $fdd_d(X)$ be the Polish space of infinite
sequences of
successive finite-dimensional subspaces in $Fin_{\Q}(X)$,
equipped with the product of the discrete
topology on $Fin_{\Q}(X)$.
The set of finite sequences of successive finite-dimensional subspaces in
$Fin_{\Q}(X)$ will be denoted by $fdd_d^{<\omega}(X)$. $\tilde{F}$
will denote a
finite sequence of successive finite-dimensional spaces,
and $(\tilde{F}_n)_{n \in \N}$ an infinite sequence of such finite sequences.
The usual notation about concatenation of finite sequences will be
used. For $S \in fdd_d(X)$,
$[S]$ will denote the linear span of $S$.

For $E,F$ finite-dimensional spaces in $X$, define the distance
$d(E,F)$  between $E$
and $F$ as the classical Hausdorff distance between the unit spheres of $E$ and
$F$.
Let $\Delta=(\delta_n)_{n \in \N}>0$. Let $A$ be a subset of $fdd_d(X)$. The
$\Delta$-expansion $A_{\Delta}$ of $A$ is the set of  sequences of successive
finite dimensional spaces $(F_k)_{k \in \N} \in fdd_d(X)$ such that
there exists $(E_k)_{k \in \N}$ in $A$ with $d(E_k,F_k) \leq \delta_k$ for
all $k \in \N$.
The following theorem was essentially proved in \cite{FR2}.

\begin{theo} \label{theo19}
Let $X$ be a Banach space with a basis. If $A$ is comeager in $fdd_d(X)$, then
for any $\Delta>0$, there exists a  sequence $(\tilde{F}_n)_{n \in \N} \in
(fdd_d^{<\omega}(X))^{\omega}$ of successive finite sequences of
successive finite dimensional
subspaces, such that all
elements of $fdd_d(X)$  passing
through infinitely many $\tilde{F}_n$'s are in $A_{\Delta}$.
\end{theo}

\pf
The proof is verbatim the same as in the case of block-sequences in
\cite{FR2} (this corresponds to Proposition \ref{comeager} in this article),
replacing blocks in
$Q(X)$ by finite-dimensional spaces in $Fin_{\Q}(X)$, and block-sequences in
$bb_d(X)$ by sequences of successive finite-dimensional subspaces in
$fdd_d(X)$.
\pff

We shall use this theorem when $A$ is in fact a property of $[x_n]_{n
\in \N}$, in that
case, each sequence $\tilde{F}_n$ can be chosen to be of length $1$,
and the formulation
becomes a little bit more tractable. It follows:

\begin{theo} \label{prop20} Let $X$ be a Banach space with an unconditional
   basis. Let $P$ be an isomorphic property of Banach spaces.
Assume that the set $\{(F_n)_{n \in \N} \in fdd(X): [F_n]_{n \in \N}
\ has \ P\}$ has the Baire Property,
and that
  $P$
is stable by
   passing to complemented subspaces and by squaring. If every sequence
in $fdd(X)$ has a subsequence whose closed linear span
satisfies $P$, then all subspaces with a successive FDD in $X$  satisfy $P$.
  \end{theo}

\pf
Let $A=\{(F_n)_{n \in \N} \in fdd(X): [F_n]_{n \in \N} \ has \ P\}$.
For small enough
$\Delta>0$,
$A_{\Delta}=A$ and
$(A^C)_{\Delta}=A^C$. In Theorem \ref{theo19}, for sets corresponding
to isomorphic properties
(such as $A$ or $A^C$), the sequence $\tilde{F}_n$
may be chosen to be of length
$1$ for each $n \in \N$. It follows from our hypotheses about $P$ that $A$
cannot be meager. For
$\tilde{E}=(E_1,\ldots,E_p) \in fdd_d^{<\omega}(X)$, denote by
$N(\tilde{E})$ the set of sequences $(F_n)_{n \in \N} \in fdd_d(X)$ such that
$F_n=E_n$ for all $n \leq p$.
As $A$ has the Baire Property, it is comeager in some open set
$N(\tilde{E})$, and without loss of generality $\tilde{E}$ is a length $1$
sequence
$(E_1)$. We now prove that $A$ is
comeager in $fdd_d(X)$.

Otherwise, $A$ is meager in some open set $N(\tilde{F})$, $\tilde{F}
\in fdd_d^{<\omega}(X)$, and
without loss of generality
$\tilde{F}$ is a length $1$ sequence $(F_1)$. Now we may find $E_2$
and $F_2$ in $Fin_{\Q}(X)$, with
$E_1<E_2$, $\dim E_2=\dim F_1$, $F_1<F_2$, $\dim F_2=\dim E_1$, and
$\max(supp(E_2))=\max(supp(F_2))$.
Let $f$ be the canonical bijection between $N((E_1,E_2))$ and
$N((F_1,F_2))$, defined by
$f((E_1,E_2)^{\frown}S)=(F_1,F_2)^{\frown}S$ for all $S \in fdd_d(X)$. It
is routine to check that $f$ is
an homeomorphism, and that for all $S \in N((E_1,E_2))$, $[f(S)]$ is
isomorphic to $[S]$; in
particular
$S
\in A$ if and only if
$f(S)
\in A$. It follows a contradiction with the fact that $A$ is meager
in $N((F_1,F_2))$ and comeager in
$N((E_1,E_2))$.

As $A$ is comeager, Theorem \ref{theo19} applies. By properties of
$P$, and because the basis of $X$
is assumed unconditional,
$A$ is stable by taking subsequences and by concatenation of disjoint
sequences.
By the same method as in the end of the proof of Proposition
\ref{theo3}, it follows that $A=fdd_d(X)$. \pff

\

One of the most important still open questions in Banach space theory is to
know whether any complemented subspace of a Banach space with an unconditional
basis must have an unconditional basis. A positive answer to this would have
many consequences, for example concerning the
Schroeder-Bernstein Property for Banach spaces (see e.g. \cite{C} for a survey).
The following corollary gives a direction for solving this question by the negative
(here we use that "spanning a subspace with an unconditional basis" is
analytic and thus has the Baire Property in $fdd_d(X)$).

\begin{coro}\label{contrexemple}  Let $X$ be a Banach space with an
unconditional basis.
Assume:

(1) every sequence in $fdd(X)$ has a subsequence which spans a
  subspace with  an unconditional basis,

(2) there exists a sequence in $fdd(X)$ which spans a
   subspace without an unconditional basis.

Then there exists a subspace $F=\oplus_{n \in \N} F_n$ of $X$ with a
   successive FDD on the basis, which has an unconditional basis, and a
   subsequence
$(G_k)_{k \in \N}$ of $(F_n)_{n \in \N}$ such that $G=\oplus_{k \in \N}G_k$,
   though complemented in $F$, does not have an unconditional basis.
\end{coro}

We conclude by discussing some of the properties that a Banach space
$X$ with (1) and (2)
must have, supposing it to exist.

Recall that a
  Banach space $X$ is said to have Gordon-Lewis {\em l.u.st.} if there
is a constant $C \geq 1$
  such that
for every finite dimensional subspace $E$ of $X$, there exists
a finite dimensional space $F$ with a $1$-unconditional basis, and maps
$T:E \rightarrow F$, $U:F \rightarrow X$, such that
$UT(x)=x$ for all $x \in E$ and such that $\norm{T}\norm{U} \leq C$.
We note that having l.u.st. is an analytic property of Banach spaces, which
is stable by passing to complemented subspaces and squaring.
As (1) implies that every sequence in $fdd_d(X)$ has a
subsequence which spans a subspace with l.u.st., it follows from Theorem
\ref{prop20} that if $X$ satisfies (1) and (2), then
every subspace of $X$ with a successive FDD must
have l.u.st..

By the Theorem of Komorowski and Tomczak-Jaegermann \cite{KT}, it
follows that $X$
must be $l_2$-saturated. Also by \cite{CK1} Theorem 3.8, every subspace of $X$
with a uniform FDD on the basis must have an unconditional basis.

Another interesting fact is that the unconditional basis for $(\oplus F_n)_{n
   \in \N}$ in the conclusion of Corollary \ref{contrexemple}  cannot be
obtained in the obvious way, that is by constructing in
each $F_n$ a $C$-unconditional basis, and proving that the sequence which is
the reunion of each basis is a $K(C)$-unconditional basis for $(\oplus F_n)_{n
   \in \N}$, for some constant $K(C)$.
In that case, any subspace $(\oplus G_k)_{k \in \N}$ associated to a
   subsequence $(G_k)_{k \in \N}$ of $(F_n)_{n \in \N}$ would inherit an
   unconditional basis (which is just a subsequence of the unconditional basis
of $(\oplus F_n)_{n \in \N}$).

\

A natural candidate for $X$ is the Orlicz sequence space $l_F$ considered by
P. Casazza and N.J. Kalton in
\cite{CK1}. It is reflexive,
has cotype $2$ and type $2-\epsilon$ for any $\epsilon>0$, and is
$l_2$-saturated. Among other interesting
properties, every subspace of $l_F$ with a uniform UFDD has an
unconditional basis. We do not know whether $l_F$ satisfies the hypotheses of Corollary
\ref{contrexemple}.

\end{document}